\documentclass[12pt,a4paper,psamsfonts]{amsart}
\usepackage{amssymb,amscd,amsxtra,calc}
\usepackage{cmmib57}

\theoremstyle{plain}
    \newtheorem{thm}{Theorem}[section]
    
    \newtheorem{claim}[thm]{Claim}
    \newtheorem{corollary}[thm]{Corollary}
    \newtheorem{lemma}[thm]{Lemma}
    \newtheorem{proposition}[thm]{Proposition}
    \newtheorem{question}[thm]{Question}
    \newtheorem{theorem}[thm]{Theorem}

\theoremstyle{definition}

    \newtheorem{remark}[thm]{Remark}

    \newtheorem{_example}[thm]{Example}
\newenvironment{example}{\begin{_example}\rm}{\end{_example}}

\theoremstyle{remark}

    \newtheorem{setup}[thm]{}

\newcommand{\PP}{\mathbb{P}}

\newcommand{\Q}{\mathbb{Q}}

\newcommand{\Z}{\mathbb{Z}}

\newcommand{\Aut}{\operatorname{Aut}}

\newcommand{\id}{\operatorname{id}}

\newcommand{\rank}{\operatorname{rank}}

\newcommand{\topol}{\operatorname{top}}

\begin{document}

\title[K3 surfaces with order 11 automorphisms]{
K3 surfaces with order 11 automorphisms}

\author{Keiji Oguiso}
\address[Keiji Oguiso]{%
\textsc{Department of Mathematics} \endgraf
\textsc{Osaka University, Toyonaka,
Osaka 560-0043, Japan, and} \endgraf
\textsc{Korea Institute for Advanced Study} \endgraf
\textsc{Hoegiro 87, Seoul, 130-722, Korea
}}
\email{oguiso@math.sci.osaka-u.ac.jp}
\author{De-Qi Zhang}
\address[De-Qi Zhang]{%
\textsc{Department of Mathematics} \endgraf
\textsc{National University of Singapore,
10 Lower Kent Ridge Road,
Singapore 119076}}
\email{matzdq@nus.edu.sg}

\begin{abstract}
In the present paper we describe the K3 surfaces admitting order 11
automorphisms and apply this to classify log Enriques surfaces
of global index 11.
$$
\text{\it This paper is
dedicated to the memory of Eckart Viehweg.}$$
\end{abstract}



\maketitle

\section{Introduction}

The purpose of this paper is to describe the family of
complex K3 surfaces with automorphisms of order 11
and apply this to classify log Enriques surfaces
of global canonical index 11 (see \cite{Z} for the definition).
We note that any automorphism of order 11 of a K3 surface
is necessarily non-symplectic, that is, acts on the space of
the global two forms non-trivially \cite{Ni}.

Throughout this paper, we consider a pair $(X,G)$
consisting of a complex projective
K3 surface $X$ and a finite group  $G$  of automorphisms on  $X$
which fits in the exact sequence:
$$
1 \rightarrow G_{N} \rightarrow G \overset{\rho}{\rightarrow}
\mu_{11n} = \langle \zeta_{11n} \rangle \rightarrow 1,
$$
where the last map $\rho$ is the natural representaion of $G$
on the space $H^{2,0}(X) = \Bbb C \omega_{X}$ and
$n$ is some positive integer. It is known that
$n \leq 6$ (\cite{Ni}, \cite{Ko1}; see also \cite{MO}).
We fix an element $g \in G$ with $\rho (g) = \zeta_{11}$, i.e.,
$$g^{*}\omega_{X} = \zeta_{11}\omega_{X}$$ and set
$$M = H^{2}(X, \Bbb Z)^{g} .$$

For simplicity of description, we also
assume that $G$ is maximal in the sense that if $(X,G')$ also
satisfies the same conditon as above for some $n'$
and $G \subseteq G'$ then $G = G'$.

In order to state our main Theorem, we first construct three types of
examples of such pairs.  We denote by $U$ and $U(m)$
the lattices defined respectively by the Gram matrix
$$\begin{pmatrix}
0 & 1 \\ 1 & 0
\end{pmatrix}, \hskip 1pc
\begin{pmatrix}
0 & m \\ m & 0
\end{pmatrix} .$$
Denote by $A_{*}$, $D_{*}$, $E_{*}$ the negative definite lattices
given by the Dynkin diagrams of the indicated types.

\begin{example}\label{Example 1} (\cite{Ko1}, \cite{MO})
Let $S_{66}$ be the K3 surface given by the Weierstrass equation
$y^{2} = x^{3} + (t^{11}-1)$, and $\sigma_{66}$ the automorphism
of $S_{66}$ given by
$$\sigma_{66}^{*}(x,y,t) =
(\zeta_{66}^{22}x, \zeta_{66}^{33}y, \zeta_{66}^{12}t) .$$
Then the pair $(S_{66}, \langle \sigma_{66} \rangle)$ gives an example of
$(X,G)$ with $n = 6$  and $G_{N} = \{1\}$, i.e., $G \simeq \mu_{66}$
and with $M \simeq U$.
\end{example}

\begin{example}\label{Example 2} (\cite{Ko1})
Consider
the rational, fibered threefold
$\varphi : \mathcal X \rightarrow \Bbb C$ defined by
$$y^{2} = x^{3} + x + (t^{11}-s)$$ and its order 22 automrphism
$\sigma$ given by
$$\sigma^{*}(x,y,t,s) = (x, -y, \zeta_{11}t, s)$$
where  $s$  is the coordinate of  $\Bbb C$. Then
$\varphi$ is a morphism smooth over $s \not= \pm \sqrt{-4/27}$
and
$${\mathcal X}_{\sqrt{-4/27}} \,\,\, (\, \simeq {\mathcal X}_{-\sqrt{-4/27}}\, )$$
has a unique singular point of type $A_{10}$.

The pair  $(\mathcal X_{0}, \langle \sigma_{44} \rangle)$,
where
$$\sigma_{44}^{*}(x,y,t) =
(\zeta_{44}^{22}x, \zeta_{44}^{11}y, \zeta_{44}^{34}t)$$
gives an example of $(X,G)$ with $n = 4$  and  $G_{N} = \{1\}$,
i.e., $G \simeq \mu_{44}$ and with $M = U$.
(The minimal resolution of)
$(\mathcal X_{s}, \langle \sigma \rangle)$  with  $s \not=0$  gives
an example of $(X, G)$ with $n = 2$ (and $G_{N} = \{1\}$), i.e.,
$G \simeq \mu_{22}$ and with $M = U$ (resp. $U \oplus A_{10}$)
if
$$s \not= 0, \pm \sqrt{-4/27} \,\,\,  \text{(\, resp. if} \,\,\, s = \pm \sqrt{-4/27} \,)$$
(cf. Remark \ref{rEx2} and the proof of Claim \ref{Claim 6} below
for the calculation of  $M$  and  $G$).
\end{example}

The following remark will help to verify the calculation
of  $G$  and  $M$  in Examples \ref{Example 1} and \ref{Example 2} above.

\begin{remark}\label{rEx2}

(1) Let  $(X, G)$  be any of the pairs in Examples \ref{Example 1} and \ref{Example 2} above
and let  $g$  be the unique order 11 element
in  $G$  satisfying  $g^* \omega_X = \zeta_{11} \omega_X$.
The natural  $G$-stable (hence  $g$-stable)
Jacobian elliptic fibration  $f : X \rightarrow {\bold P}^1$,
with  $t$  as the inhomogeneous coordinate of the base space,
is the only  $g$-stable elliptic fibration on  $X$
(cf. the first paragraph in the proof of Proposition \ref{Proposition 3} below.)

(2) The fixed locus (point wise)
$X^{g}$  is equal to the union of a smooth rational curve
in the type  $I_{11}$  fiber  $X_{t = 0}$  and two points
on the type  $II$  fiber  $X_{t = \infty}$ (resp. the union of
the smooth fiber  $X_{t = 0}$  and two points on the type $II$
fiber  $X_{t=\infty}$), when  $X$  is equal to  ${\mathcal X}_{\sqrt{-4/27}}$
(resp. any of other cases in Examples \ref{Example 1} and \ref{Example 2}).

(3) For any  $s \ne 0, \pm \sqrt{-4/27}$, four
surfaces
$$S_{66}, \,\, {\mathcal X}_0,  \,\, {\mathcal X}_{\sqrt{-4/27}}, \,\, {\mathcal X}_s$$
are not isomorphic to one another.
\end{remark}

\begin{example}\label{Example 3}

Let us consider the following
three series of rational Jacobian elliptic surfaces:

(1) $j^{(1)} : J^{(1)} \rightarrow \Bbb P^{1}$,
defined by the Weierstrass equation
$$y^{2} = x^{3} + (t-1) $$
whose singular fibers are $J^{(1)}_{1}$ of Kodaira's type $II$
and $J^{(1)}_{\infty}$ of Kodaira's type $II^{*}$;

(2) $j^{(2)} : J^{(2)} \rightarrow \Bbb P^{1}$,
defined by the Weierstrass equation
$$y^{2} = x^{3} + x + (t-s)$$
with $s \not= \pm \sqrt{-4/27}$,
whose singular fibers are $J^{(2)}_{\alpha}$, $J^{(2)}_{\beta}$
of Kodaira's type $I_1$,
and $J^{(2)}_{\infty}$ of Kodaira's type $II^{*}$,
where  $t = \alpha, \beta$  are two distinct non-zero roots
of the discriminant  $\Delta(t) = 4 + 27(t-s)^2$; and

(3) $j^{(3)} : J^{(3)} \rightarrow \Bbb P^{1}$,
defined by the Weierstrass equation
$$y^{2} = x^{3} + x + (t-s)$$
with  $s = \sqrt{-4/27}$
whose singular fibers are $J^{(3)}_{0}$, $J^{(3)}_{2s}$
of Kodaira's type $I_{1}$,
and $J^{(3)}_{\infty}$ of Kodaira's type $II^{*}$.

Let $p^{(i, e)} : P^{(i, e)} \rightarrow \Bbb P^{1}$ be a
non-trivial principal homogeneous space of
$j^{(i)} : J^{(i)} \rightarrow \Bbb P^{1}$ given
by an element $e$ of order 11 in  $(J^{(i)})_{0}$.
(For the basic results on the principal homogeneous space
of rational Jacobian elliptic fibrations, see \cite[Chapter V, Section 4]{CD}.)
Then $p^{(i, e)} : P^{(i, e)} \rightarrow \Bbb P^{1}$
is a rational elliptic surface with a multiple fiber
of multiplicity $11$ over $0$ (of type $I_{0}$ in the cases $i = 1, 2$
and of type $I_{1}$ in the case $i = 3$).

Let $Z^{(i, e)}$ be the log Enriques surface of index $11$
obtained by the composite of the blow up
at the intersection of the components
of multiplicities 5 and 6 in $(P^{(i, e)})_{\infty}$,
which is of Kodaira's type $II^{*}$, and the blow down of
the proper transform of $(P^{(i, e)})_{\infty}$. Let $X^{(i,e)}$
be the global canonical cover of $Z^{(i, e)}$ and $G^{(i,e)}$
the Galois group of this covering. Then, each of these pairs
$(X^{(i,e)}, G^{(i,e)})$ gives an example of $(X,G)$ with
$n = 1$  and $G_{N} = \{1\}$, i.e.,
$G \simeq \mu_{11}$ and with $M = U(11)$ (see Lemma \ref{Lemma 9} below
to verify the calculation of  $G$  and  $M$).
\end{example}

Our main result is as follows:

\begin{theorem}[Main Theorem]\label{ThA}

Under the notation above, the following are true.

\begin{itemize}
\item[(1)]
We have  $G_{N} = \{1\}$  so that  $G \simeq \mu_{11n}$
and  $g$  is unique and of order  $11$.

\item[(2)]
$M$ is isomorphic to either one of $U$, $U \oplus A_{10}$
or $U(11)$.

\item[(3)]
In the case where $M \simeq U$ or $U \oplus A_{10}$,
$(X,G)$ is isomorphic to either $(S_{66}, \langle \sigma_{66} \rangle)$,
$(\mathcal X_{0}, \langle \sigma_{44} \rangle)$, or
$(\mathcal X_{s}, \langle \sigma \rangle)$ ($s \not=0$) in Examples $\ref{Example 1}$ and $\ref{Example 2}$.

Moreover, $M \simeq U \oplus A_{10}$
if and only if $(X,\langle g \rangle)$ is isomorphic to
$$(\mathcal X_{\sqrt{-4/27}},  \langle \sigma^{2} \rangle) \,\,\,
(\, \simeq (\mathcal X_{-\sqrt{-4/27}}, \langle \sigma^{2} \rangle) \,).$$

\item[(4)]
In the case where $M \simeq U(11)$,
$(X,G)$ is isomorphic to one of  $(X^{(i,e)}, G^{(i,e)})$
in Example $\ref{Example 3}$.
\end{itemize}

\end{theorem}

Combining the main Theorem \ref{ThA} with Remark \ref{rEx2},
we obtain the following, where a log Enriques surface
is maximal if, by definition, any birational morphism
$Z' \rightarrow Z$  from another log Enriques surface $Z'$
must be an isomorphism.

\begin{corollary}

Maximal log Enriques surfaces of global index $11$
are isomorphic to either a  $Z^{(i, e)}$ in Example $\ref{Example 3}$
or $\overline{\mathcal X}_{\sqrt{-4/27}}/\langle g \rangle$,
where  $\overline{\mathcal X}_{\sqrt{-4/27}}$  is the surface
obtained from the surface  ${\mathcal X}_{\sqrt{-4/27}}$
in Example $\ref{Example 2}$ with the unique  $g$-fixed curve contracted.

\end{corollary}

\begin{remark}
(1) In the main Theorem \ref{ThA} $(3)$ and $(4)$ and Examples \ref{Example 2} and \ref{Example 3},
the pairs  $(X,G)$  parametrized by $s$ and $-s$, are isomorphic
to each other.
In particular, the pair $(X,G)$  with $M \simeq U \oplus A_{10}$
is unique up to isomorphisms.

(2) By the main Theorem \ref{ThA}, the pairs $(X,G)$  are not finitely many
any more and move in a 1-dimensional
(non-isotrivial) family, which is one of
the main difference from the previous works \cite{MO}, \cite{OZ2}, \cite{Xi}, \cite{Ko1} \cite{Ko2}
concerning larger non-symplectic group actions.
Indeed, calculating the $J-$invariant and combining with the fact
that the pair $(X,G)$ with $\text{ord}(G) = 40$ and its
elliptic fiber space structure are both unique \cite{MO}, we find that
the family $\varphi : \mathcal X \rightarrow \Bbb C$
given in Example \ref{Example 2} is not isotrivial. Similarly, the uniqueness
of the Jacobian elliptic fiber space structure
on a rational surface
shows that the family given in Example \ref{Example 3} is also not isotrivial.

(3) One can also explain the reason why  $(X,G)$'s form a 1-dimensional
family from the view point of the period mapping.
Since for generic $(X,G)$, the transcendental lattice $T_{X}$ is
of rank 20 and isomorphic to either $U^{2} \oplus E_{8}^{2}$
or $U \oplus U(11) \oplus E_{8}^{2}$; further, the eigenspace
with respect to the eigenvalue $\zeta_{11}$ of
the action  $g$  on  $T_X \otimes {\Bbb C}$
in which the period $\Bbb C\omega_{X}$ should lie is
two dimensional.
Conversely a generic one dimensional subspace in this eigen space
gives periods of K3 surfaces with order 11 automorphisms
$g$ by the surjectivity of the period mapping \cite{BPV}.

(4) In our classification, we make use of the invariant part
$M$  of the  $g$-action on  $H^2(X, {\Bbb Z})$,
instead of the Neron Severi lattice $S_{X}$ which always contains $M$
and certainly equals  $M$ if  $X$ is generic in the family.
However, for special  $X$, $S_{X}$ is probably larger than $M$.
So, in our classification, the determination of the Neron Severi
lattice \cite{Sh}, which is one of the hardest and most important problems
concerning algebraic surfaces,
remains unsettled.  The reason why we describe
the result according to $M$ rather than $S_{X}$
is that on the one hand,
the Neron Severi lattices are quite unstable under deformations,
for instance, in the case of the family of quartic K3 surfaces,
and on the other hand, it turns out
that the invariant part $M$ is fairly stable under deformation
at least in our case.
\end{remark}

A group $G$ is called a {\it $K3$ group}, if $G \le \Aut(X)$
for some complex $K3$ surface $X$.

\begin{proposition}\label{Mon}
No sporadic finite simple group which is different from the Monster group ${\bf M}$,
contains all finite $K3$ groups as its subgroups.
\end{proposition}

\begin{question}\label{Q}
Can we embed every finite $K3$ group into the Monster simple group {\bf M}?
\end{question}

\begin{remark}
After this work was done and motivated by Mukai's embedding of
{\it all} finite symplectic $K3$ groups into the sporadic simple Mathieu group $M_{23}$ ($\le {\bf M}$)
and the observation in Proposition \ref{Mon}, the above
Question \ref{Q} crossed our minds.
We planned to solve this question and include the current paper as part of the new project
\cite{OZ3}.
However, this project is unexpectedly complicated and we have not yet completed it.
So we decide to publish the current paper as an independent paper.

After the current paper was written in 1999,
there have been much progress, especially in positive characteristic,
among which is the very significant work of Dolgachev-Keum \cite{DK} where the authors
successfully extended Mukai's classification of finite symplectic $K3$ groups
to positive characteristics. See also \cite{Z2} for a partial survey.
\end{remark}

{\bf Acknowledgement.}  This paper is finalized during the
second author's visits to Japan in June 1998 and March 1999,
and appeared as  {\it arXiv:math/9907020}.
The first named author is supported by JSPS Program 22340009 and
KIAS Scholar Program, and the second named author
is supported by an ARF of NUS.
The authors are very thankful to the JSPS-NUS programme
for the financial support, and the referee for pointing out an error and typos.

\section{Proof of the main Theorem}

We now prove the main Theorem \ref{ThA}.
We employ the same notation introduced in \S 1 freely.
Let  $N$  be the orthogonal lattice of $M$ in $H^{2}(X, \Bbb Z)$.
Then $N$ is $g-$stable and $T_{X} \subseteq N$ and $M \subseteq S_{X}$.
For a lattice $L$, we denote by $L^{*}$ the dual (over) lattice
Hom$(L, \Bbb Z)$. For a positive integer $I$,
we denote by $\varphi(I)$ the cardinality of
the multiplicative group  $(\Bbb Z/I)^{\times}$.

\begin{lemma}\label{Lemma 1}
We have  $G_{N} = \{1\}$. In particular, $g$ is unique and is of
order $11$.
\end{lemma}

\begin{proof}

Suppose the contrary that $|G_N| \ge 2$.
Since $\varphi(11) \, \vert \, \text{rank} \, T_{X}$,
$\text{rank} \, T_{X}$ is either $10$ or $20$.
According \cite[the list]{Xi2} and its notation,
$$\rank \, T_X = \rank T_Y  \le 22 - (c + 1)$$
with $c \ge 8$ (resp. $c \ge 12$) when $|G_N| \ge 2$ (resp. $|G_N| \ge 3$).
Thus $|G_N| = 2$ and $\rank \, T_X = 10$. Hence $\rank S_X = 12$.

Write $G_{N} = \langle \iota \rangle \simeq \Bbb Z/2$
Then $|G| = 22$.
Since $G_N \lhd G$, we have $G = \langle h \rangle \simeq \Bbb Z/22$.
We may assume that $g = h^2$ and $\iota = h^{11}$.
By the topological Lefschetz fixed point formula,
we have the diagonalization $\iota^{*} \vert (S_{X} \otimes {\Bbb C})
= \text{diag}[I_4, -I_8]$, relative to some basis.
By considering minimal polynomial of $\zeta_{11}$ over $\Q$,
we have either
$g^* | S_X \otimes {\Bbb C} = I_{12}$ (identity matrix), or
$$g^* | S_X \otimes {\Bbb C} = \text{diag}[\zeta_{11}, \zeta_{11}^2, \dots, \zeta_{11}^{10}, 1, 1] .$$

If the second case for $g^*$ occurs, then simultaneously diagonalize
$g^*$ and $\iota^*$ on $S_X \otimes {\Bbb C}$, we would get a diagonalization
of $(g \circ \iota)^*$ whose diagonal entries consist of a few $\pm 1$ and between $6$ and
$8$ entries of $22nd$ primitive roots of the unity, which is impossible
because $g \circ \iota$ is of order $22$ and the Euler number $\varphi(22) = 10$.

If the first case for $g^*$ occurs, then we get the following diagonalizations,
relative to two possibly different bases (up to re-ordering):
$$(g \circ \iota)^{*} \vert (S_{X} \otimes {\Bbb C})
= \text{diag}[I_4, -I_8], \,\,\,
(g \circ \iota)^{*} \vert (T_{X} \otimes {\Bbb C})
= \text{diag}[\zeta_{11}, \zeta_{11}^2, \dots, \zeta_{11}^{10}] .$$
Thus $\chi_{topol}(X^{g \circ \iota}) = -3$ by the
topological Lefschetz fixed point formula. In particular,
$X^{g \circ \iota}$ contains a curve.  On the other hand,
since $(g \circ \iota)^{11} = \iota$, we have
$X^{g \circ \iota} \subseteq X^{\iota}$, so that
$X^{g \circ \iota}$ consists of finitely many points,
a contradiction.

\end{proof}

\begin{lemma}\label{Lemma 2}

$M$ is isomorphic to either $U$, $U(11)$ or
$U \oplus A_{10}$.

\end{lemma}

\begin{proof}

Since  $M$ is a primitive sublattice of the unimodular
lattice  $H^{2}(X, \Bbb Z)$, we have a natural isomorphism
$M^{*}/M \simeq N^{*}/N$.  Noting that  $N^{g^{*}} = \{0\}$, we
can apply the same argument
as in \cite[Lemmas (1.1), (3.2)]{MO} and \cite[Lemmas (1.2), (1.3)]{OZ2}
for the pair $(M, N)$ (instead of $(S_{X}, T_{X})$  there) to get
$\varphi(11)=10 \, \vert \, \text{rank} \, N$ and
$$M^{*}/M \simeq N^{*}/N \simeq (\Bbb Z/11)^{\oplus s}$$ for some
integer $s$ with $0 \leq s \leq \text{rank} \, N/10$. Since
$\text{rank} \, N \leq 21$, we have
$$(\text{\rm rank}(M), \text{\rm det}(M)) = (22- \text{\rm rank}(N),
-11^s) =$$
$$(2, -1), (2, -11), (2, -11^2), (12, -1), (12, -11).$$

By the classification of indefinite unimodular even lattices,
the case  $(\text{\rm rank}(M)$, $\text{\rm det}(M))$
$= (12, -1)$  is impossible
and in the case  $(\text{\rm rank}(M), \text{\rm det}(M)) = (2,-1)$
we have  $M = U$.

By \cite{RS}, a $p$-elementary ($p > 2$)
even hyperbolic lattice of rank $> 2$,
is determined uniquely by its rank and discriminant.
So, $M = U \oplus A_{10}$  when
$(\text{\rm rank}(M), \text{\rm det}(M)) =
(12, -11)$.

Suppose that  rank $M = 2$.  Write  $M = (a_{ij})$,  where
$a_{11} = 2a, a_{22} = 2c, a_{12} = a_{21} = b$
for integers  $a, b, c$.
Then det$M \equiv 0, -1$ (mod 4)  and hence
the case  $(\text{\rm rank}(M), \text{\rm det}(M)) = (2, -11)$
is impossible.  We consider the case where
$(\text{\rm rank}(M)$, $\text{\rm det}(M)) = (2, -11^2)$.
Note that  $M^*$  is generated by a  ${\Bbb Z}$-basis
$$(\varepsilon_1  \,\, \varepsilon_2) = (e_1 \,\, e_2) M^{-1}
= (e_1 \,\, e_2) (-1/11^2)(b_{ij})$$  where
$$b_{11} = 2c, b_{22} = 2a, b_{12} = b_{21} = -b .$$
Here  $e_i$'s  form the basis of  $M$  with  $(a_{ij})$
as the intersection matrix.  Since  $M^*/M = (\Bbb Z/11)^{\oplus s}$,
each  $b_{ij}$ (and hence each  $a_{ij}$) is divisible by  $11$.
So  $M = M_1(11)$  with an indefinite even unimodular lattice  $M_1$.
Thus  $M = U(11)$  under a suitable basis.

\end{proof}

\begin{proposition}\label{Proposition 3}

Assume that $M \simeq U$. Then $(X,G)$ is isomorphic
to either $(S_{66}, \langle \sigma_{66} \rangle)$,
$(\mathcal X_{0}, \langle \sigma_{44} \rangle)$, or
$(\mathcal X_{s}, \langle \sigma \rangle)$ ($s \not= 0, \pm \sqrt{-4/27}$)
in Examples $\ref{Example 1}$ and $\ref{Example 2}$.

\end{proposition}

\begin{proof}

If $n \geq 3$, the result follows from \cite[Main Theorem]{MO}.
Let us consider the case $n \leq 2$.
Since $M \simeq U$, $X$ admits a $g-$stable Jacobian
fibration  $f : X \rightarrow \Bbb P^{1}$ by \cite{OS}.
Let $E$ and $C$ be a general fiber of $f$ and
the unique $g-$stable section
of $f$. Here the uniqueness of the $g-$stable section follows from
the fact that if $C'$ is also a  $g-$stable section then
$[C'] = a[C] +b[E]$ and
$$((aC+bE).E) = 1, \,\,\, (aC+bE)^{2} = -2 .$$
We see then these equalities imply $a = 1$ and $b = 0$.

Let $\overline{g}$ be the automorphism of the base space
$\Bbb P^{1}$ induced by $g$. Since there are no elliptic curves
admitting Lie automorphism of order 11, $\overline{g}$ is also
of order 11. We may then adjust an inhomogeneous coordinate $t$
of ${\Bbb P}^{1}$   so that
$({\Bbb P}^{1})^{\overline{g}} = \{0, \infty\}$.
We note that $X_{0}$ and $X_{\infty}$ are both irreducible,
because the irreducible component $R$ of $X_{0}$ meeting $C$
is $g-$stable so that $\text{rank} \, M \geq 3$ unless $R = X_{0}$.

Since $g^{*}\omega_{X} = \zeta_{11}\omega_{X}$, an easy
local coordinate calculation shows that
neither of  $X_0, X_{\infty}$  is of Kodaira's type $I_{1}$.
Moreover, noting that  $g$ permutes the other singular fibers, we have
$$24 = \chi_{topol}(X) = \chi_{topol}(X_{0}) + \chi_{topol}(X_{\infty}) + 11m$$
for some positive integer $m$.
Thus after suitable change of inhomogeneous
coordinate $t$ if necessary, $(X_{0}, X_{\infty})$ is
of type  $(I_{0}, II)$  and the set of the
other singular fibers is either

(1) $\{X_{\zeta_{11}^{i}} | 0 \leq i \leq 10\}$,
all of Kodaira's type $II$,

(2) $\{X_{\zeta_{11}^{i}}, X_{\zeta_{11}^{j} \alpha} |
0 \leq i, j \leq 10\}$ ($\alpha \notin \mu_{11}$), all
of Kodaira's type $I_{1}$, or

(3) $\{X_{\zeta_{11}^{i}} | 0 \leq i \leq 10\}$, all of
Kodaira's type $I_{2}$.

\end{proof}

\begin{claim}\label{Claim 4}

The case $(3)$ can not happen.

\end{claim}

\begin{proof}

Assuming the contrary that Case(3) occurs,
we denote by  $R$  the irreducible component
of $X_{1}$ meeting $C$. Since
$$S := \sum_{0 \leq i \leq 10}g^{i}(R)$$ is $g-$stable, we have
$[S] = a[C] + b[F]$.  Now  $(S.F) = 0$  implies that  $a = 0$
and hence  $S = b[F]$. This leads to
$$-22 = (S)^{2} = (bF)^{2} = 0$$
which is a contradiction.

\end{proof}

\begin{claim}\label{Claim 5}

The case $(1)$ can not happen under the assumption that $n \leq 2$.

\end{claim}

\begin{proof}
Assuming the contrary that Case (1) happens, we
will determine the Weierstrass equation
$$y^{2} = x^{3} + a(t)x + b(t)$$ of $f : X \rightarrow \Bbb P^{1}$.
Since the singular fibers of $f$ are all of type $II$,
the $J-$function
$$J(t) := 4a(t)^{3}/(4a(t)^{3}+27b(t)^{2}) = 0$$
as a rational function.
Thus, $a(t) = 0$ and the equation is $y^{2} = x^{3} + b(t)$.

Let us consider the discriminant divisor
$$\Delta(t) = 27b(t)^{2} .$$
Since the singular fibers of $f$ over $t \not= \infty$ are
$X_{\zeta_{11}^{i}}$ and these are all of type $II$, we have
$\Delta(t) = c(t^{11} - 1)^{2}$ for some nonzero constant $c$.
Then $b(t) = c'(t^{11}-1)$ for some nonzero constant $c'$.
Changing  $x,y$  by suitable multiples, we finally
find that $f$ is given by the equation
$$y^{2} = x^{3} + (t^{11}-1)$$
which is isomorphic to $S_{66}$ in Example \ref{Example 1}. In particular,
$G \simeq \mu_{66}$ by \cite{MO}. Thus $n = 6$, a contradiction.
The referee pointed out that the argument above is similar to
\cite[(5.1)]{Ko1}; we keep this argument for readers' convenience.

\end{proof}

\begin{claim}\label{Claim 6}

Assume that $f : X \rightarrow \Bbb P^{1}$
satisfies the condition of the case $(2)$
and $M \simeq U$ and $n \leq 2$. Then $f : X \rightarrow \Bbb P^{1}$
is isomorphic to a Jacobian elliptic fibration given by a Weierstrass
equation
$$y^{2} = x^{3} + x + (t^{11}-s)$$  for some
$s \not= 0, \pm \sqrt{-4/27}$, and under this isomorphism,
we have  $G \simeq \langle \sigma \rangle$,
where
$$\sigma^{*}(x,y,t) = (x, -y, \zeta_{11}t) .$$  In particular, $n =2$.

\end{claim}

\begin{proof}

Again we will determine the Weierstrass equation
$$y^{2} = x^{3} + a(t)x + b(t)$$ of $f : X \rightarrow \Bbb P^{1}$,
where  $a(t), b(t)$  are polynomials in  $t$.
First note that $\text{deg} \, a(t) \leq 8$ and $\text{deg} \, b(t) \leq 12$
by the canonical bundle formula.
Since  $f$ has singular fibers
$$\{X_{\zeta_{11}^{i}}, X_{\zeta_{11}^{j} \alpha} \, | \, 0 \leq i, j \leq 10 \}$$
of type $I_{1}$, the discriminant divisor
$\Delta(t)$  is equal to
$$\delta(t^{11}-\alpha^{11})(t^{11}-1)$$
for some non-zero constant  $\delta$.
Since the $J-$function
$$J(t) = 4a(t)^3/\Delta(t)$$
is  $\overline{g}-$invariant,
$a(t)$ (and hence $b(t)$) are also $\overline{g}-$semi invariant.
Thus
$$a(t) = At^m, \,\, b(t) = t^n(B_1+B_2t^{11})$$ where
$A, B_i$  are constants, $m \le 8, n \le 12$, and  $n \le 1$
when  $B_2 \ne 0$.  Comparing coefficients of the equality
$$\Delta(t) = 4a(t)^3 + 27 b(t)^2$$ we see that
$$a(t) = A, \, b(t) = B_1 + B_2t^{11} .$$
Noting that  $A \ne 0$  because of the existence of
singular fibers of type  $I_1$. We have also
$B_2 \ne 0$, otherwise, $X$  is birational to a product of a fibre
and the parameter space $\PP^1$ and hence is not a $K3$ surface, absurd!
We can, by a suitable coordinate change,
normalize the Weierstrass equation of  $X$  as
$$X = {\mathcal X}_s : \,\, y^{2} = x^{3} + x + (t^{11} - s).$$
Here  $s$  is a constant, and  $s \ne 0$  for otherwise  $n = 4$ by \cite{MO}.

Conversely, by the standard algorithm
to finding out the singular fibers \cite{Ne}, we see that this elliptic surface
$\mathcal X_{s}$ has 22 singular fibers of type $I_{1}$ and
a singular fiber of type $II$ if and only if $s \not= \pm \sqrt{-4/27}$.
Moreover, $\mathcal X_{s}$ admits an automorphism $g_{s}$ of order $11$
given by
$$g_{s}^{*}(x,y,t) = (x, y, \zeta^{11} t) .$$
Since $g$ and $g_{s}$ make the fibration  $f$  and the section
$C$  stable and satisfy
$$g^{*}\omega_{X} = g_{s}^{*}\omega_{X} = \zeta_{11}\omega_{X}$$ we have
$g = g_{s}$.  Now the condition that  $n \le 2$  implies
that  $n = 2$  and  $G \simeq \mu_{22}$,
by the maximality of  $G$  and by
considering  $$G_s = \langle  g_s, \iota_s  \rangle \simeq \mu_{22}$$
where  $$\iota_{s}^{*}(x,y,t) = (x,-y,t)$$  acts on $f$ as the
involution around  $C$.

Write  $G = \langle  g, \iota  \rangle$  with an involution  $\iota$.
Since $\iota \circ g = g \circ \iota$, we see that
$C$ and $f$ are both $\iota-$stable, and
$0$  and  $\infty$  are two  $\iota$-fixed points.
In other words, $\iota$ does not switch $0$ and $\infty$,
because the fibres $X_0$ and $X_{\infty}$ are of different types: $I_0$, $II$.
If $\iota$ acts
on the base space $\Bbb P^{1}$ as an involution, $G$  permutes the
22 singular fibers of type $I_{1}$  as well as the 22 roots
of the discriminant divisor
$$\Delta(t) = 4 + 27(t^{11} - s)^2$$
whence $s = 0$, a contradiction.
Thus, $\iota$ is the involution of $f$ around $C$, i.e.,
$\iota = \iota_s$.  This means $G = G_{s}$  and we are done.
\end{proof}

This completes the proof of Proposition \ref{Proposition 3}.

\par \vskip 1pc
Next we consider the case where $M \simeq U \oplus A_{10}$.
In this case, $M = S_{X}$ and $\text{rank} \, T_{X} = \varphi(11) = 10$.
So  $(X, G)$  is equivariantly isomorphic to the pair
$({\mathcal X}_{\sqrt{-4/27}}, \langle  \sigma  \rangle)$
in Example \ref{Example 2}, by \cite[Theorem 2]{OZ2} and by making use of
the maximality of  $G$ as in the previous paragraph.
This also proves the main Theorem \ref{ThA} in the case of
$M = U$  or  $U \oplus A_{10}$.

\par \vskip 1pc
Finally we consider the case where $M \simeq U(11)$. As before,
since $U(11)$ represents zero, $X$ admits a  $g-$stable elliptic fibration
$f : X \rightarrow \Bbb P^{1}$ and the induced action $\overline{g}$
on the base space is of order 11. We adjust an inhomogeneous coordinate $t$
of the base so that  $(\Bbb P^{1})^{\overline{g}} = \{0, \infty\}$. We need
further coordinate change later, but we always keep this condition.

\begin{lemma}\label{Lemma 7} After a suitable coordinate change,
$f$ satisfies either one of the following three cases.
\begin{itemize}
\item[(1)]
$X_{0}$ is smooth and $g\vert X_{0}$ is a translation of order $11$;
the remaining singular fibers are $X_{\zeta_{11}^{k}}$
$(0 \leq k \leq 10)$ and these are all of type $II$.

\item[(2)]
$X_{0}$ is smooth and $g\vert X_{0}$ is a translation of order $11$;
the remaining singular fibers are $X_{\zeta_{11}^{k}}$ and
$X_{\zeta_{11}^{k}}\alpha$ $(0 \leq k \leq 10$ and $\alpha \notin
\mu_{11})$ and these are all of type $I_{1}$.

\item[(3)]
$X_{0}$ is of Type $I_{11}$ and $g\vert X_{0}$ is a translation of order $11$
(which permutes the fiber components cyclically);
the remaining singular fibers are $X_{\zeta_{11}^{k}}$
and these are all of type $I_{1}$.
\end{itemize}

Moreover, in all three cases, $X_{\infty}$ is of type $II$
with  $X^g = (X_{\infty})^{g} = \{P_{1}, P_{2}\}$, where $P_{1}$ is the
singular point of $X_{\infty}$.  The action of $g$ around $P_{i}$ is of
type $1/11(5,7)$ if $i=1$ and $1/11(2,10)$ if $i=2$.

\end{lemma}

\begin{proof} The proof is almost identical to the situation
where  $M \supseteq U$, except that $f$ does not
admit $g-$stable sections and we use the assumption that
$M \simeq U(11)$ and the fact that $X^{g}$ is smooth. The type of
the action is determined by an elementary local coordinate calculation
of the normalization of $X_{\infty}$ and the fact that
$g^{*}\omega_{X} = \zeta_{11}\omega_{X}$.
Actually, we have one more possible case in which $X_{0}$ is smooth with
$g\vert X_{0} = \id$  and  $X_{\infty}$ is of type $II$ with
$(X_{\infty})^{g} = \{P_{1}, P_{2}\}$.
But then the relatively minimal
model of $X/\langle g \rangle \rightarrow \Bbb P^{1}/\langle \overline{g} \rangle$ is a rational
elliptic surface with no multiple fibers and hence has a section $C$.
Now the pullback on  $X$  of  $C$  is a  $g$-stable section, which
contradicts $M \simeq U(11)$.
\end{proof}

Note that the fibration $f$  on  $X$  induces an elliptic fibration
$f' : X/\langle g \rangle \rightarrow
\Bbb P^{1}/\langle \overline{g} \rangle$  on the
quotient surface, a log Enriques surface of index 11.
Let $S \rightarrow X/\langle g \rangle$ be the
minimal resolution.
Then the proper transform $D_0$ of $X_{\infty}/\langle g \rangle$
is a $(-1)$-curve on  $S$. This is because that the total transform $D$ of
$X_{\infty}/\langle g \rangle$ is a non-relatively minimal fibre of an elliptic fibration on the smooth surface $S$;
to be precise,
every irreducible component $D_i$ ($\ne D_0$) of this fibre $D$ is a curve with
self-intersection $\le -2$, and at least one curve say $D_1$ has $D_1^2 \le -3$ since
$P_1$ (and also $P_2$) is not a rational double point.

We let $c : S \rightarrow T$ be the contraction of this  $(-1)$ curve and
$\overline{f} : T \rightarrow \Bbb P^{1}$ the induced
relatively minimal rational
elliptic fibration.  We immediately get the following
lemma from the construction.

\begin{lemma}\label{Lemma 8}
According to the cases $(1), (2), (3)$ in Lemma $\ref{Lemma 7}$, the singular fibers
of $\overline{f}$ are:

\begin{itemize}
\item[(1)]
$T_{0}$ of type ${}_{11}I_{0}$, $T_{\infty}$ of type $II^{*}$, and
$T_{1}$ of type $II$.

\item[(2)]
$T_{0}$ of type ${}_{11}I_{0}$, $T_{\infty}$ of type $II^{*}$, and
$T_{1}$ and $T_{\alpha^{11}}$ of type $I_1$.

\item[(3)]
$T_{0}$ of type ${}_{11}I_{1}$, $T_{\infty}$ of type $II^{*}$, and
$T_{1}$ of type $I_1$.
\end{itemize}
\end{lemma}

Note that we can recover  $(X, f)$  in Lemma \ref{Lemma 7} easily
from  $(T, \overline{f})$  in Lemma \ref{Lemma 8}.  Indeed,
let $\overline{f} : T \rightarrow \Bbb P^{1}$ be a
relatively minimal rational elliptic surface with one of the properties
(1), (2) and (3) in Lemma \ref{Lemma 8}. Blow up the point of the intersection
of the components of multiplicities 5 and 6 in $T_{\infty}$ and then
contract the two connected components of
the proper transform of $T_{\infty}$.
We now get a rational elliptic surface $f' : S \rightarrow \Bbb P^{1}$
with two singular points of types $1/11(5,7)$ and $1/11(2,10)$ and with
$11K_{S}$  linearly equivalent to $0$.  Let $X \rightarrow S$ be the
global canonical  ${\Bbb Z}/11$-cover of $S$.  Then
$(X, f)$,  where  $f$   is induced from  $\overline{f}$,
fits corresponding cases in Lemma \ref{Lemma 7}.

Moreover, if we let  $F$ be the Galois group  Gal$(X/S)$,
then we have:

\begin{lemma}\label{Lemma 9}
$(X, F)$ satisfies the condition in the second paragraph
of the Introduction. Further, $M \simeq U(11)$ and  $F$ is maximal.
In particular, in the situation of Lemma $\ref{Lemma 7}$, one has
$(X, F) = (X, \langle  g  \rangle)$.
\end{lemma}

\begin{proof}
The first assertion is clear. We use the notation, like $g$, $M$ in the Introduction for $F$.

If $M \not \simeq U(11)$, then $M$ is either $U$
or $U \oplus A_{10}$ by Lemma \ref{Lemma 2}. However, then $X^{g}$ contains a curve
by the main Theorem \ref{ThA} (proved already when  $M \supseteq U$)
and Remark \ref{rEx2},
which contradicts the fact that the canonical covering map is
\'etale in codimension one. Thus $M \simeq U(11)$.

Next, we show that $F$ is maximal. Assume that $F \subset H$
and $H$ also satisfies the condition in the Introduction.
By Lemma \ref{Lemma 1}, $H_N = \{1\}$.
By \cite{MO}, it is enough to eliminate the case where $H = \langle h \rangle \simeq \Z/22$.

Assume the contrary that this case happens.
We may assume that the order $11$ element $g := h^2$ is as in the Introduction.
Since
rank $M = 2$, $F_N = \{1\}$ and
$\text{rank} \,  T_{X}$ is
either $10$ or $20$, we have either
$$h^{*}\vert T_{X} \otimes \Bbb C =
\text{diag}[-\zeta_{11}^j| 1 \le j \le 10]$$  and
$h^{*}\vert S_{X} \otimes \Bbb C$  equals one of:
$$\text{diag}[1, \pm 1, \zeta_{11}^{j} | 1 \le j \le 10], \,\,\,\,
\text{diag}[1, \pm 1, -\zeta_{11}^{j} | 1 \le j \le 10]$$
in the case where $\text{rank} \, T_{X} = 10$, or
$$h^{*}\vert T_{X} \otimes \Bbb C =
\text{diag}[-\zeta_{11}^{j} | 1 \le j \le 10]^{\oplus 2}, \,\,\,\,
h^{*}\vert S_{X} \otimes \Bbb C =
\text{diag}[1, \pm 1]$$ in the case where $\text{rank} \, T_{X} = 20$.

Since $h(X^{g}) = X^{g}$ and $X^{h} \subseteq X^{g}$,
we have $X^{h} = \{P_{1}, P_{2}\}$, noting
that the actions of $g$ around two points  $P_i$
are different.  Thus the topological Lefschetz formula
shows that the only possible case is:
$$h^{*}\vert T_{X} \otimes \Bbb C =
\text{diag}[-\zeta_{11}^j | 1 \le j \le 10], \,\,\,\,
h^{*}\vert S_{X} \otimes \Bbb C =
\text{diag}[1, -1, \zeta_{11}^{j} |1 \le j \le 10] .$$

Let  $\iota := h^{11}$. Then,
$$\iota^{*}\vert T_{X} \otimes \Bbb C = -I_{10}, \,\,\,\,
\iota^{*}\vert S_{X} \otimes \Bbb C = \text{diag}[I_{11}, -1] .$$
In particular, $\chi_{\topol}(X^{\iota}) = 2$.
This, together with the fact that
$\iota^{*}\omega_{X} = -\omega_{X}$, implies that
$X^{\iota}$ consists of smooth curves and
at least one of them is a smooth rational curve, say $C$.
Write the (disjoint) irreducible decomposition of $X^{\iota}$  as
$$X^{\iota} = C \cup E_{1} \cup ... \cup E_{m} .$$
Since $g \circ \iota = \iota \circ g$, the $g$ acts on the set
$\{C, E_{1},..., E_{m}\}$.

First assume that $g(C) \not= C$.
Then $g^{i}(C)$ would be mutually disjoint 11 rational curves with
$${\Bbb Q}\langle g^{i}(C) \rangle \subseteq
S_{X}^{\iota^{*}} \otimes {\Bbb Q}$$
where both sides of the inclusion are of rank 11, whence
they are equal.  However,
$S_{X}^{\iota^{*}}$  then contains no ample classes, a contradiction.
Thus $g(C) = C$ and $P_{1}, P_{2} \in C$. But,
this can not happen, because
the action of $g$ around $P_{i}$ are of types $1/11(5,7)$ if $i=1$ and
$1/11(2,10)$  if  $i = 2$, and there are no  $a \in \{5,7\}, b \in \{2,10\}$
with $a + b \equiv 0 \, (\text{mod} \, 11)$.
Therefore, $F$  is maximal and Lemma \ref{Lemma 9} is proved.
\end{proof}

Now the only remaining task is to describe rational elliptic surfaces
with the property (1), (2), or (3) in Lemma \ref{Lemma 8}.
However, each of these is obtained as a
principal homogeneous space of a Jacobian rational elliptic surface
$j : J \rightarrow \Bbb P^{1}$  whose singular fiber type
is equal to one of the three types in Example \ref{Example 3}.
Now a similar (and easier) calculation shows that
the Weierstrass equation of  $j : J \rightarrow \Bbb P^{1}$
is the same as one of those in Example \ref{Example 3}.
This completes the proof of the main Theorem \ref{ThA}.

\begin{setup}
{\bf Proof of Proposition \ref{Mon}}
\end{setup}
Let $I$ be either $54$ or
a prime number $\le 19$, then the order-$I$ cyclic group
$\mu_I$ acts purely non-symplectically on
some $K3$ surface and hence is a $K3$ group (cf. \cite[Main Theorem 3]{MO}).
Among the $26$ sporadic simple groups in \cite{Atlas}, only the Monster ${\bf M}$ contains
all such $\mu_I$ as subgroups (neither the baby Monster $B$ nor the
Mathieu group $M_{23}$ contains $\mu_{54}$ as its subgroup).
This proves the proposition.


\begin{thebibliography}{99}

\bibitem{BPV}
W. Barth, C. Peters and A. Van de Ven, Compact complex surfaces,
Springer-Verlag, 1984.

\bibitem{Atlas} J. H. Conway, R. T. Curtis, S. P. Norton, R. A. Parker and R. A. Wilson,
\textit{Atlas of finite groups.
Maximal subgroups and ordinary characters for simple groups. With computational assistance from J. G. Thackray},
Oxford University Press, Eynsham, 1985.

\bibitem{CD}
F. R. Cossec and I. V. Dolgachev, Enriques surfaces I, Progress
in Mathematics, Vol. {\bf 76}, Birkhauser 1989.

\bibitem{DK}
I. V. Dolgachev and J. Keum,
Finite groups of symplectic automorphisms of $K3$ surfaces in positive characteristic,
Ann. of Math. (2) {\bf 169} (2009), no. 1, 269--313.

\bibitem{Ko1}
S. Kondo, Automorphisms of algebraic  $K3$  surfaces which
act trivially on Picard groups, J. Math. Soc. Japan, {\bf 44} (1992), 76-98.

\bibitem{Ko2}
S. Kondo,
The maximum order of finite groups of automorphisms of $K3$ surfaces,
Amer. J. Math. {\bf 121} (1999), no. 6, 1245--1252.

\bibitem{MO}
N. Machida and K. Oguiso, On  $K3$  surfaces admitting finite
non-symplectic group actions, J. Math. Sci. Univ. Tokyo, {\bf 5} (1998),
273-297.

\bibitem{Ne}
A. Neron, Mod\'eles minimaux des vari\'et\'es ab\'eliennes
sur les corp locaux et globaux,
Publ. Math. I.H.E.S. {\bf 21} (1964).

\bibitem{Ni}
V. V. Nikulin, Finite automorphism groups of Kahler  $K3$  surfaces,
Trans. Moscow Math. Soc. {\bf 38} (1980), 71-135.

\bibitem{OS}
K. Oguiso and J. Sakurai,
Calabi-Yau threefolds of quotient type, Asian J. Math. {\bf 5} (2001), no. 1, 43--77.

\bibitem{OZ1}
K. Oguiso and D. -Q. Zhang, On the most algebraic  $K3$  surfaces
and the most extremal log Enriques surfaces,
Amer. J. Math. {\bf 118} (1996), 1277-1297.

\bibitem{OZ2}
K. Oguiso and D. -Q. Zhang,
On Vorontsov's theorem on $K3$ surfaces with non-symplectic group actions,
Proc. Amer. Math. Soc. {\bf 128} (2000), no. 6, 1571--1580.

\bibitem{OZ3}
K. Oguiso and D. -Q. Zhang, Finite automorphism groups of
$K3$  surfaces, in preparation.

\bibitem{RS}
A. N. Rudakov and I. R. Shafarevich, Surfaces of type K3 over
fields of finite characteristics, Current Problems in Math. {\bf 18}
(1981), 115-207.

\bibitem{Sh}
T. Shioda, On the Mordell-Weil lattices,
Comment. Math. Univ. Sancti Pauli, {\bf 39} (1990), 211-240.

\bibitem{Xi}
G. Xiao, Non-symplectic involutions of a  $K3$  surface, preprint 1996.

\bibitem{Xi2}
G. Xiao, Galois covers between $K3$ surfaces, Ann. Inst. Fourier, Grenoble,
{\bf 46} (1996), 73-88.

\bibitem{Z}
D. -Q. Zhang, Logarithmic Enriques surfaces, I, II;
J. Math. Kyoto Univ. {\bf 31} (1991), 419-466; {\bf 33} (1993), 357-397.

\bibitem{Z2}
D. -Q. Zhang,
Automorphisms of K3 surfaces,
Proceedings of the International Conference on Complex Geometry and Related Fields,
pp. 379--392, AMS/IP Studies in Advanced Mathematics, 39, Amer. Math. Soc., Providence, RI, 2007.

\end{thebibliography}
\end{document}